\LoadClass[11pt]{article}

\documentclass[11pt]{article}

\usepackage{amsmath, amssymb, amsthm} 
\usepackage{listings,verbatim,float }
\usepackage{graphicx}
\usepackage{tikz}
\usepackage{geometry}
\usepackage{forest}
\usepackage{hyperref}  
\usepackage{ragged2e} 
\usepackage{fancyhdr}

\usepackage{setspace} 
\usepackage[style=authoryear,backend=biber]{biblatex}
\title{Quantifying the Balinski-Young Theorem:\\ Structure and Probability of Quota Violations in Divisor Methods for Three States}
\author{Tyler C. Wunder and Joseph W. Cutrone}

\newcommand{\thetitle}{}
\newcommand{\theauthor}{}
\renewcommand{\title}[1]{\renewcommand{\thetitle}{#1}}

\renewcommand{\author}[1]{\renewcommand{\theauthor}{#1}}
\AtBeginDocument{%
 \begin{center}
 {\LARGE\textbf{\thetitle}\par}
 \vskip 1.7em
 {\large\theauthor}
 \vskip 2em
 \end{center}
}

\newtheorem{theorem}{Theorem}[section]
\newtheorem{corollary}{Corollary}[theorem]
\newtheorem{lemma}[theorem]{Lemma}

\newtheorem{deflem}[theorem]{Definition/Lemma}
\newtheorem{definition}[theorem]{Definition}

\newtheorem{remark}[theorem]{Remark}

\newenvironment{contact}{
 \vspace{0.5em}
 \begin{center}
 \begin{tabular}{@{}l}
}{
 \end{tabular}
 \end{center}
 \vspace{0.5em}
}

\newcommand{\email}[1]{Email:~\url{#1}}

\newcommand{\ffo}{\lfloor q_1 \rfloor}
\newcommand{\fft}{\lfloor q_2 \rfloor}

\newtheorem{innercustomthm}{Theorem}
\newenvironment{customthm}[1]
 {\renewcommand\theinnercustomthm{#1}\innercustomthm}
 {\endinnercustomthm}

 \hypersetup{colorlinks=true,
 linkcolor=blue,
 citecolor=blue,
 urlcolor=blue,
 pdfauthor={\theauthor},
 pdftitle={\thetitle}}

\renewenvironment{abstract}{
 \begin{center}
 \small\textbf{Abstract}\\[8pt]
 \begin{minipage}{.85\textwidth}
}{
 \end{minipage}
 \end{center}
 \vspace{1em}
}

\newcommand{\noi}{\noindent}

\newcommand{\LL}{\left}
\newcommand{\RR}{\right}

\title{Quantifying the Balinski-Young Theorem:\\ Structure and Probability of Quota Violations in Divisor Methods for Three States}
\author{Tyler C. Wunder and Joseph W. Cutrone}

\begin{document}

\frenchspacing
\begin{abstract}
The apportionment problem asks how to assign representation to states based on their populations. That is, given census data and a fixed number of seats, how many seats should each state be assigned? Various algorithms exist to solve the apportionment problem, but by the Balinski-Young Impossibility Theorem, every such algorithm will be flawed in some way. This paper focuses on divisor methods of apportionment, where the possible flaws are known as quota violations.\\
\\
This paper presents a detailed analysis of quota violations that arise under divisor methods for three states. By focusing on the three-state case, our analysis makes the consequences of the Balinski-Young theorem particularly transparent and enables a precise classification of quota violations. The study focuses on quota violations in the Adams, Jefferson, Dean, and the Huntington-Hill methods when allocating $M$ seats, but extends to a wider class of divisor functions. Theoretical results are proved about the behavior of these methods, particularly focusing on the types of quota violations that may occur, their frequency, and their geometry. The key results of the paper are algebraic tests to detect quota violations which are employed to understand the geometry of violations and construct a probability function which calculates the likelihood of such violations occurring given an initial three state population vector whose components follow varying distributions. 

\end{abstract}


\section{Introduction}
The \textit{apportionment problem} - dividing a fixed number of indivisible resources among states in proportion to their populations - lies at the intersection of political science, combinatorics, and probability. A key class of apportionment methods are divisor methods, which are a group of simple methods that avoid unfair results known as apportionment paradoxes \parencite{end:b:meeting-ideal}. Divisor methods have been widely employed in the United States for assigning seats in House of Representatives, with the Jefferson, Webster, or Huntington-Hill method used for most of the country's history. The Huntington-Hill method is currently used to apportion the House of Representatives. \parencite{FYS_Textbook}

While divisor methods avoid paradoxes, these methods do not always satisfy the \textit{quota rule}, which requires that each state receive an allocation equal to either the floor or ceiling of its exact proportional share. This phenomenon of imperfect apportionment is in fact guaranteed by the Balinski-Young Impossibility Theorem \parencite{end:b:meeting-ideal}, which establishes that no apportionment method can simultaneously satisfy the quota rule and avoid paradoxes.

Since divisor methods can violate quota, and as flawed methods are unavoidable, natural questions to ask are: What is the underlying structure of such violations? How often will a quota violation happen? And, can quota violations be predicted in advance given states' populations? This paper answers these questions for three states. Specifically, we investigate the frequency, structure, and geometry of quota violations under divisor methods in three states in the Adams, Jefferson, Webster, Dean, and Huntington-Hill methods. The paper characterizes how quota violations occur and provides formulas to calculate the probability that a divisor method violates quota. Overall, this paper
provides a detailed analysis of quota violations in the simpler case of three states, which ultimately highlights and corroborates the Balinski Young Impossibility Theorem.

This paper builds on related work about probabilities and analysis of apportionment within social choice theory and the mathematics of elections, such as similar probabilistic work on simplexes concerning the Jagiellonian compromise found in \parencite{zyczkowski2004votingeuropeanunionsquare} and \parencite{SlomczynskiZyczkowski2007HO}; work on systemic bias such as in \parencite{SchusterPukelsheimDrtonDraper2003}, \parencite{SchwingenschloglDrton2004}, \parencite{DrtonSchwingenschlogl2005}, and \parencite{Ichimori2012}; and work on the probabilities of voting paradoxes such as in \parencite{gehrlein2017elections} and \parencite{pandit-cutrone-fairness-2025}.

\section{Preliminary Notation and Definitions}

\begin{definition}\label{Divisor_def} (Divisor Method) 
 Given a strictly increasing divisor function $d: \mathbb{Z}_{\geq 0} \to \mathbb{R}_{\geq 0}$, populations $(p_1, \dots, p_n)$ and a fixed number of seats $M$, an apportionment function $A(p_1, \dots, p_n) = (a_1, \dots, a_n)$, is then defined as follows: 
\begin{enumerate}
\item[\emph{1.}] \emph{Assign each state zero seats.}
\item[\emph{2.}] \emph{For each state calculate their \emph{priority value}, $ \frac{p_i}{d(r_i)}$, where $r_i$ is the number of seats currently assigned to state $i$.} 
\item[\emph{3.}] \emph{Assign to the state with the highest priority value another seat.}
\item[\emph{4.}] \emph{Repeat steps (2.) and (3.) until $M$ seats have been apportioned.}
\end{enumerate}

If $d(0)=0$, use the convention that division by zero results in infinity, and as such all states immediately receive one seat. Equivalently, if $d(0) = 0$, replace the first step with assigning each state one seat rather than zero. 
\end{definition}
The most common and arguably reasonable divisor methods are the ``five workable methods'' \parencite{b-young-methods,end:b:meeting-ideal, Huntington-theory} defined by the following divisor functions: 
\begin{center}
\begin{tabular}{|l|l|}
  \hline
  Name & Divisor Function \\ \hline
  Adams & $d(s) = s$ \\
  Jefferson & $d(s) = s+1$\\
  Webster & $d(s) = s+0.5$\\
  Huntington-Hill & $d(s) = \sqrt{s(s+1)}$ \\
  Dean & $d(s) = \dfrac{2s(s+1)}{2s+1}$ \\ \hline
  
\end{tabular}
\end{center}

Important values are the \emph{standard quotas} and \emph{standard divisors}. Let $P = \sum_{i=1}^n p_i$ be the total population. Then, the \emph{standard divisor} is defined as $SD = \frac{P}{ M }$, and state $i$'s \emph{standard quota} is $q_i = \frac{p_i}{SD} .$ The standard quota is a percentage $\frac{p_i}{P}$ of $M$ and is a state's theoretical fair share of the number of seats. As such each state should be apportioned roughly their standard quota. Since the standard quota is rarely an integer, states should be apportioned either their \emph{upper quota}, $\lceil q_i \rceil$, or their \emph{lower quota}, $\lfloor q_i \rfloor$. 

A \emph{quota violation} occurs whenever a state is assigned more seats than their upper quota or fewer seats than their lower quota. Specifically, a \emph{lower quota violation} occurs when $a_i < \lfloor q_i \rfloor$ and an \emph{upper quota violation} occurs when $a_i > \lceil q_i \rceil$. Informally, a quota violation is when a state is apportioned too many or too few seats.

The Balinski-Young Theorem, seen in \parencite{end:b:meeting-ideal}, guarantees that all apportionment methods either have quota violations or have strange, unfair occurrences known as \emph{apportionment paradoxes}. Divisor methods have quota violations but avoid paradoxes \parencite{end:b:meeting-ideal}.

\section{Results}
The following is a summary of the results of the paper:
\begin{customthm}{\ref{no-up}} For $M$ seats and $n=3$ states, if $A$ is a divisor apportionment method with divisor function $d(s) \in [s, s+\frac{1}{2}]$ for all $s$, then upper quota violations cannot occur. In particular, upper quota violations do not occur in the Adams, Dean, and Huntington-Hill methods for $n = 3$ states. \end{customthm}

\begin{customthm}{\ref{extra}} (\textit{Classification of Lower Quota Violations in 3 States})
Order the populations $p_1< p_2<p_3$. If $A$ is the Adams, Dean, or the Huntington Hill method and $A(p_1,p_2,p_3)$ has a quota violation, the final apportionment must be $(a_1,a_2,a_3) =(\lceil q_1 \rceil, \lceil q_2 \rceil, \lfloor q_3 \rfloor - 1) $, where $q_i$ is the standard quota of state $i$. That is, the quota violation must be a lower quota violation occurring in the third state.\end{customthm}

\begin{customthm}{\ref{q.v_criteria_test}} \textit{(Lower Quota Violation Criteria Test)} Let $A$ be the Adams, Dean, or the Huntington-Hill method, with divisor function $d(s)$. Order the populations such that $p_1 <p_2 < p_3$ and calculate the standard quotas $q_i$. Then the apportionment $A(p_1, p_2, p_3)$ on $M$ seats has a lower quota violation if and only if all three of the following statements are true:
\begin{enumerate}
\item $q_3 d(\lfloor q_1 \rfloor) < q_1 d( \lfloor q_3 \rfloor - 1)$
\item $q_3 d(\lfloor q_2 \rfloor) < q_2 d( \lfloor q_3 \rfloor - 1)$
\item $\lceil q_1 \rceil + \lceil q_2 \rceil + \lfloor q_3 \rfloor = M+1 $ (equivalently for $q_1, q_2 \notin \mathbb{Z}$, $\lfloor q_1 \rfloor + \lfloor q_2 \rfloor + \lfloor q_3 \rfloor = M-1) $.
\end{enumerate}

\end{customthm}

\begin{customthm}{\ref{lower_prob}} \textit{(Probability of Lower Quota Violation)} 
 For $M$ seats and $n=3$ states, let $q_1$ and $q_2$ be the standard quotas of states 1 and 2 with joint pdf $f$, and define $q_3=M-q_1-q_2$. For $q_1, q_2 < q_3$ and for apportionment method $A$ equal to Adams, Dean, or Huntington-Hill, the probability of a lower quota violation is given by:
\[ \sum_{\ffo,\fft} \iint_{\bigtriangleup_l(\ffo,\fft)} f(\ffo+x,\fft+y)\,dx\,dy \]
where $\bigtriangleup_l(\ffo,\fft)$ denotes the feasible region of the Lower Quota Violation Criteria Test (\ref{q.v_criteria_test}) defined in Theorem \ref{l_tri_def}.\end{customthm}

\begin{customthm}{\ref{extra2}}\textit{(Classification of Upper Quota Violations in 3 States)} Order the populations $p_1 < p_2<p_3$. If $A$ is the Jefferson method and $A(p_1, p_2, p_3)$ has a quota violation, the final apportionment must be $(a_1,a_2,a_3)= (\lfloor q_1 \rfloor, \lfloor q_2 \rfloor, \lceil q_3 \rceil +1)$, where $q_i$ is the standard quota of state $i$. That is, the quota violation must be an upper quota violation occurring in the third state. \end{customthm}

\begin{customthm}{\ref{uq.v_criteria_test}} \textit{(Upper Quota Violation Criteria Test)}
Let $A$ be the Jefferson method, with divisor function $d(s)$. Order the populations such that $p_1 <p_2 < p_3$ and calculate the standard quotas $q_i$. Then the apportionment $A(p_1, p_2, p_3)$ on $M$ seats has an upper quota violation if and only if all three of the following statements are true:
\begin{enumerate}
\item $q_3 d( \lfloor q_1 \rfloor) > q_1 d(\lceil q_3 \rceil)$
\item $q_3 d( \lfloor q_2 \rfloor) > q_2 d(\lceil q_3 \rceil)$
\item $\lfloor q_1 \rfloor + \lfloor q_2 \rfloor + \lceil q_3 \rceil = M-1$ (equivalently for $q_3 \notin \mathbb{Z}$, $\ffo + \fft + \lfloor q_3 \rfloor = M-2$.)
\end{enumerate}
\end{customthm}

\begin{customthm}{\ref{up_prob}} (Probability of Upper Quota Violation)
 For $M$ seats and $n=3$ states, let $q_1$ and $q_2$ be the standard quotas of states 1 and 2 with joint pdf $f$, and define $q_3=M-q_1-q_2$. For $q_1, q_2 < q_3$ and the Jefferson method, the probability of an upper quota violation is given by:
 
\[ \sum_{\ffo,\fft} \iint_{\bigtriangleup_u(\ffo,\fft)} f(\ffo+x,\fft+y)\,dx\,dy \]
where $\bigtriangleup_u(\ffo,\fft)$ denotes the feasible region of the Upper Quota Violation Criteria Test (\ref{uq.v_criteria_test}) defined in Theorem \ref{u_tri_def}.
\end{customthm}
 

\section{Background Information}

The following are simple lemmas, known results, and common assumptions that are utilized in this paper. Throughout the paper, let $M$ be the number of seats and $n$ the number of states. Let $A(p_1,\ldots,p_n)$ be the divisor method apportionment for state populations $p_i$, with $A$ an arbitrary divisor method. 

\subsection{Assumptions}
Assume there are no ties in priority values; in particular this includes assuming no ties in states' populations. Additionally, to guarantee at least one seat per state assume $M \geq n$.

\subsection{Introductory Lemmas and Definitions}
\begin{lemma}\label{neutrality} \parencite{end:b:meeting-ideal} (Symmetry) A divisor method is preserved under reordering. That is, if $A(p_1, \dots , p_n)=(a_1, \dots, a_n)$ then for any re-indexing, $(1, \dots, n) \mapsto (\sigma(1), \dots, \sigma(n))$, A($p_{\sigma(1)}, \dots , p_{\sigma(n)}) = (a_{\sigma(1)}, \dots, a_{\sigma(n)})$. 
\end{lemma}

\begin{lemma}\label{proportionality} \parencite{end:b:meeting-ideal} (Homogeneity) A divisor method is preserved under proportionality. If $A(p_1, \dots , p_n)=(a_1, \dots, a_n)$ then for any scalar, $\lambda \in \mathbb{R}_{>0}$, $A(\lambda p_1, \dots , \lambda p_n)=(a_1, \dots, a_n)$. In particular, apportioning the standard quotas is the same as apportioning the populations: $A(p_1, \dots, p_n) = A(q_1, \dots, q_n)$. .
\end{lemma}

\begin{deflem}
 \label{Hill=HH} Define a divisor apportionment method $\tilde{A}(p_1,\ldots,p_n)$ on $n$ states and $M$ seats by finding a modified divisor $D \in \mathbb{R}_{>0}$ such that $D$ causes $\sum_{i=0}^n a_i = M$, for $a_i$ calculated as follows:
\begin{enumerate}
\item For each state, calculate a modified quota $m_i = \frac{p_i}{D}.$
\item For each state, assign seats by modified-rounding: round up if $m_i \geq d( \lfloor m_i \rfloor )$; round down if less.
\end{enumerate}
\end{deflem}

\noi The equivalence of the apportionment methods $A(p_1,\ldots,p_n)$ in Definition \ref{Divisor_def} and $\tilde{A}(p_1,\ldots,p_n)$ is classical. See \parencite{end:b:meeting-ideal}, \parencite{Huntington-theory}, and \parencite{Pukelsheim2017}. 

\begin{theorem}\label{no_up_and_low} \parencite{end:b:meeting-ideal, Boundsmal} No divisor method will have both upper and lower quota violations for any $M$.\end{theorem}

\begin{theorem}\label{nvv} \parencite{end:b:meeting-ideal} (No Violation Theorem) A summary of cases when no quota violations occur is as follows:
\begin{itemize}
 \item[(a)] For $n=2$ states, no upper or lower quota violations occur in any divisor methods. 
 \item[(b)] For $n=3$ states, the Webster method is the unique method such that upper or lower quota violations do not occur. 
 \item[(c)] The Adams Method is the unique method where upper quota violations do not occur for all $n$. 
 \item[(d)] The Jefferson Method is the unique method where lower quota violations do not occur for all $n$. 
\end{itemize}
\end{theorem}


\section{Structure and Probability of Quota Violations for Three States in the Adams, Dean, and Huntington-Hill Methods}
We start by providing a detailed study of lower quota violations in the Adams, Dean, and Huntington-Hill methods when apportioning three states. The Webster and Jefferson methods have no lower quota violations by Theorem \ref{nvv}.

\subsection{Characterization of Quota Violations in Three States}
By the following theorem, only lower quota violations can occur in the Adams, Dean, and Huntington-Hill methods for three states.

\begin{theorem}\label{no-up} For $M$ seats and $n=3$ states, if $A$ is a divisor apportionment method with divisor function $d(s) \in [s, s+\frac{1}{2}]$ for all $s$, then upper quota violations cannot occur. In particular, upper quota violations do not occur in the Adams, Dean, and Huntington-Hill methods for $n = 3$ states. \end{theorem}
\begin{proof} 
Let $p_i$ and $q_i$ be state $i$'s populations and standard quotas respectively. Then $A(p_1,p_2,p_3) = A(q_1, q_2, q_3) = (a_1,a_2,a_3)$. Let $d(s) \in [s, s+\frac{1}{2}]$ for all $s$. Assume for contradiction there exists an upper quota violation and, without loss of generality, the violation is on state one. Then, $a_1 \geq \lceil q_1 \rceil + 1$, and if $q_1$ is not an integer $a_1 \geq \ffo + 2$.

Write each $q_i$ as $\lfloor q_i \rfloor + d_i$, where $d_i$ is the decimal part. Since an upper quota violation exists by assumption, by Theorem \ref{no_up_and_low}, no lower quota violation exists. Therefore $a_2 \geq \lfloor q_2 \rfloor$ and $a_3 \geq \lfloor q_3 \rfloor$. Then if $q_1$ is not an integer:
\begin{align*}
\ffo + \fft + \lfloor q_3 \rfloor + d_1 + d_2 + d_3 &= q_1 + q_2 + q_3 \\
&= M \\
&= a_1 + a_2 + a_3 \\
&\geq \lceil q_1 \rceil + 1 +\fft + \lfloor q_3 \rfloor \\
&= \ffo + 2 + \fft + \lfloor q_3 \rfloor
\end{align*}
Which implies $\sum d_i \geq 2$. Since $d_i \in [0,1)$, at least two $d_i$ are greater than $\frac{1}{2}$. Since $d(s) \in [s, s+\frac{1}{2}]$, there are at least two roundups when assigning apportionment for the modified divisor $D=SD$, resulting in two states assigned their upper quota and one state with an unknown assignment of either its lower or upper quota. 

Similarly, if $q_1$ is an integer, then:
\begin{align*}
\ffo + \fft + \lfloor q_3 \rfloor +d_2 + d_3 &= q_1 + q_2 + q_3 \\
&= a_1 + a_2 + a_3 \\
&\geq \lceil q_1 \rceil + 1 +\fft + \lfloor q_3 \rfloor \\
&= \ffo + 1 + \fft + \lfloor q_3 \rfloor
\end{align*}
Which implies $d_2+d_3 \geq 1$ and therefore at least one of $d_2$, $d_3$ is greater than $\frac{1}{2}$. This results in one of those two states apportioned its upper quota and state one apportioned its upper quota. The final apportionment is then two states assigned their upper quota and one state with an unknown assignment of either its lower or upper quota.

As such, for $q_1$ an integer or not, two states are assigned their upper quota and one state has an unknown assignment of either their lower or upper quota. If this third state is assigned its lower quota, pick $d=SD$ and there is no upper quota violation. If instead the upper quota is assigned, pick $D > SD$. This cannot create an upper quota violation as increasing $D$ reduces apportionments. Both cases contradict the assumption that an upper quota violation occurred. 
\end{proof}

The proof above allows for lower quota violations to occur in $n=3$ states and indeed they do. The following three lemmas completely classify the lower quota violations in the Adams, Dean and Huntington-Hill apportionment methods.

\begin{lemma}\label{violaton_shape} In a divisor apportionment method with $M$ seats, $n=3$ states, and standard quotas $q_i$, if a lower quota violation occurs then the apportionment of the state with the violation is $\lfloor q_i \rfloor - 1$ and the other two states are assigned their upper quotas. \end{lemma}

\begin{proof}
Assume without loss of generality that the lower quota violation occurs on state 3. The claim is that the final apportionment is then $a_1 = \lceil q_1 \rceil$, $a_2 = \lceil q_2 \rceil$, and $a_3 = \lfloor q_3 \rfloor -1$. Since the first two states cannot be apportioned more than the ceiling of their quota by Theorem \ref{no_up_and_low}, $a_1 \leq \lceil q_1 \rceil$ and $a_2 \leq \lceil q_2 \rceil$. Similarly state 3 cannot have an apportionment larger than $\lfloor q_3 \rfloor -1$ by the existence of a lower quota violation, so $a_3 \leq \lfloor q_3 \rfloor-1$.

Write $q_i = \lfloor q_i \rfloor + d_i$, where $d_i$ is the decimal part. Then $M = \sum q_i = \sum \lfloor q_i \rfloor + d_i = \sum \lfloor q_i \rfloor + \sum d_i$. In particular, $\sum d_i = M - \sum \lfloor q_i \rfloor$ is an integer and since each $d_i \in [0,1), \sum d_i \in \{0,1,2\}.$ 

Then $$a_1 + a_2 = M - a_3 \geq M - (\lfloor q_3 \rfloor -1) = \sum \lfloor q_i \rfloor + \sum d_i - \lfloor q_3 \rfloor +1 = \lfloor q_1 \rfloor + \lfloor q_2 \rfloor + \sum d_i + 1$$

Proceed for each case of $\sum d_i \in \{0,1,2\}:$

If $\sum d_i =2, a_1 + a_2 \geq \lfloor q_1 \rfloor + \lfloor q_2 \rfloor +3 = \lceil q_1 \rceil + \lceil q_2 \rceil +1$. This results in an upper quota violation which is not possible by Theorem \ref{no_up_and_low} since there exists a lower quota violation by assumption. 

If $\sum d_i =0$, all of the quotas are integers. However, if the quotas are all integers, this contradicts the assumption that a lower quota violation exists.

If $\sum d_i =1,a_1 + a_2 \geq \lfloor q_1 \rfloor + \lfloor q_2 \rfloor +2 \geq \lceil q_1 \rceil + \lceil q_2 \rceil$. This is the largest possible value of $a_1 + a_2$ so equality must hold. Then $a_3 = M - a_1 - a_2 = \sum \lfloor q_i \rfloor + \sum d_i - (\lfloor q_1 \rfloor + 1) - (\lfloor q_2 \rfloor +1) = \lfloor q_3 \rfloor -1$. The final apportionment is then $a_1 = \lceil q_1 \rceil$, $a_2 = \lceil q_2 \rceil$, and $a_3 = \lfloor q_3 \rfloor -1$ and, as this is the only possible case, the claim is proven. 
 
\end{proof}

In the above lemma, the state with the lower quota violation was arbitrary. The next lemma will justify calling it state three when populations are ordered $p_1 < p_2 < p_3$:

\begin{lemma}\label{vio_on_largest}
In the Adams, Webster, and Huntington-Hill methods for $n=3$ states, lower quota violations can only occur on the state with the largest population. 
\end{lemma}
\begin{proof}
Begin by ordering the states such that $q_1 < q_2 < q_3$. Assume for the sake of contradiction there is a lower quota violation on state $j$ for $j = 1$ or $2$. Then, by Lemma \ref{violaton_shape} state three is apportioned its upper quota. Thus, the priority value that assigns state three its upper quota is greater than the priority value that would put state $j$ at its lower quota: 
\[ \frac{q_j}{d(\lfloor q_j \rfloor -1)} < \frac{q_3}{d(\lfloor q_3 \rfloor)}.\]
This simplifies to
\[q_j d(\lfloor q_3 \rfloor) < q_3d(\lfloor q_j \rfloor -1) \]
Note that $q_3 > q_j > 1$ as a lower quota violation is assumed and states are guaranteed at least one seat.

The next lemma, Lemma \ref{ineq1}, shows the above inequality implies $q_j > q_3$, a contradiction.
\end{proof}

\begin{lemma}\label{ineq1} Let $d$ be the divisor function of the Adams, Dean, or Huntington-Hill method. Then for $a,b \in [1, \infty)$, $a \cdot d(\lfloor b\rfloor) < b \cdot d(\lfloor a \rfloor - 1)$ implies $a \geq b$.
\end{lemma}
\begin{proof}
For the Dean method, let $d(s) = \frac{2s(s+1)}{2s +1} $. First, notice if $a < 2$, then $a \cdot d(\lfloor b\rfloor) < b \cdot d(\lfloor a \rfloor)$ never holds. Thus $a \geq 2$. Additionally, this implies the statement is true if $b<2$. The proof for $a,b \geq 2$ will be done by showing the contrapositive. Let $2 \leq a\leq b$. Let $n = \lfloor a \rfloor, m = \lfloor b \rfloor$. Note $m \geq n \geq 1$. Define $f(x) = \dfrac{x}{d( \lfloor x \rfloor-1)}$ and $g(x) = \dfrac{x}{d(\lfloor x\rfloor)}$. Then, it suffices to show $f(a) \geq g(b) $. 

To do so, first notice for $k \in \mathbb{Z}, k \geq 1$, on each interval $[k, k+1)$, both $f(x)$ and $g(x)$ are increasing as their derivatives are positive. Since $f$ is increasing, $f(a) \ge f(n)$ on $[n, n+1)$. Similarly, since $g$ is increasing, $g(b) \leq g(m + 1^-)$ where $m+1^-$ indicates approaching $m+1$ from the left on $[m,m+1)$. Then, by simplifying and rearranging:
\[ g(m+1^{-}) = \dfrac{(m+1)(m-0.5)}{m(m-1)} = f(m+1 ). \]
As $m \geq n$, both are integers, and $\dfrac{(m+1)(m-0.5)}{m(m-1)}$ is a decreasing function for $m \in (1, \infty)$. It follows that $f(n) \geq g(m+1^-)$. Hence, $f(a) \geq f(n) \geq g(m+1^-) \geq g(b)$. \\
\\
For the Huntington-Hill method or the Adams method, the proof is similar.
 
\end{proof}

Summarizing the results of Lemma \ref{vio_on_largest} and Lemma \ref{violaton_shape} gives the following:

\begin{theorem}\label{extra}(\textit{Classification of Lower Quota Violations in 3 States})
Order the populations $p_1< p_2<p_3$. If $A$ is the Adams, Dean, or Huntington-Hill method and $A(p_1,p_2,p_3)$ has a quota violation, the final apportionment must be $(a_1,a_2,a_3) =(\lceil q_1 \rceil, \lceil q_2 \rceil, \lfloor q_3 \rfloor - 1) $, where $q_i$ is the standard quota of state $i$. That is, the quota violation must be a lower quota violation occurring in the third state. \end{theorem}

\subsection{Lower Quota Violation Criteria Test for 3 States}
Next a test is provided to check for lower quota violations on three states. This provides an analytic check for a quota violation which will then be used to find the probability of a lower quota violation for 3 states. 

\begin{theorem}\label{q.v_criteria_test} \textit{(Lower Quota Violation Criteria Test)} Let $A$ be the Adams, Dean, or Huntington-Hill method, with divisor function $d(s)$. Order the populations such that $p_1 <p_2 < p_3$ and calculate the standard quotas $q_i$. Then the apportionment $A(p_1, p_2, p_3)$ on $M$ seats has a lower quota violation if and only if all three of the following statements are true:
\begin{enumerate}
\item $q_3 d(\lfloor q_1 \rfloor) < q_1 d( \lfloor q_3 \rfloor - 1)$
\item $q_3 d(\lfloor q_2 \rfloor) < q_2 d( \lfloor q_3 \rfloor - 1)$
\item $\lceil q_1 \rceil + \lceil q_2 \rceil + \lfloor q_3 \rfloor = M+1 $ (equivalently for $q_1, q_2 \notin \mathbb{Z}$, $\lfloor q_1 \rfloor + \lfloor q_2 \rfloor + \lfloor q_3 \rfloor = M-1) $.
\end{enumerate}

\end{theorem}
\begin{proof}
First assume both $q_1, q_2 \notin \mathbb{Z}$. Next, assume $\lfloor q_i \rfloor >0$ for all $i$. Then the Lower Quota Violation Criteria Test is equivalent to: 
\begin{enumerate}
 \item $\dfrac{q_3}{d( \lfloor q_3 \rfloor - 1)} < \dfrac{q_1}{d\lfloor q_1 \rfloor )}$
 \item $\dfrac{q_3}{d( \lfloor q_3 \rfloor - 1)} < \dfrac{q_2}{d\lfloor q_2 \rfloor )}$
 
 \item $\lfloor q_1 \rfloor + \lfloor q_2 \rfloor + \lfloor q_3 \rfloor = M-1 $.
 \end{enumerate}
To prove this biconditional statement, first assume a lower quota violation from $A(p_1,p_2,p_3)$. Then, by Theorem \ref{extra}, the apportionment is $(a_1, a_2, a_3) = (\lceil q_1 \rceil, \lceil q_2 \rceil, \lfloor q_3 \rfloor - 1) $. 

Then the priority value to apportion state three its lower quota is less than the priority values to apportion states one and two their upper quota. Thus,
$$\frac{q_3}{d( \lfloor q_3 \rfloor - 1)} < \frac{q_1}{d(\lfloor q_1 \rfloor )} \text{ \,and\, } \frac{q_3}{d( \lfloor q_3 \rfloor - 1)} < \frac{q_2}{d(\lfloor q_2 \rfloor)}.$$
Finally, from the apportionment,
\[ M = a_1 + a_2 + a_3 = \lceil q_1 \rceil+ \lceil q_2 \rceil+ \lfloor q_3 \rfloor - 1 = \lfloor q_1 \rfloor + \lfloor q_2 \rfloor + \lfloor q_3 \rfloor + 1 \] and Condition (3) then follows.

Now assume (1)-(3) are true and continue to assume $\lfloor q_i \rfloor > 0$. Start to compute $A(p_1, p_2, p_3)$= $A(q_1, q_2, q_3)$ and consider the apportionment algorithm after all states are apportioned $\lfloor q_i \rfloor -1$ seats (guaranteed by by Lemma \ref{violaton_shape}). Then, by (3), $\lfloor q_1 \rfloor + \lfloor q_2 \rfloor + \lfloor q_3 \rfloor = M-1$ and there are 4 seats left to apportion.
By (1) and (2), the priority values of states 1 and 2 at their ceilings are higher than the priority value of state 3 at its floor. Thus, the remaining four seats are apportioned to states 1 and 2 and the final apportionment is $(a_1, a_2, a_3) = (\lceil q_1 \rceil, \lceil q_2 \rceil, \lfloor q_3 \rfloor - 1)$, resulting in a lower quota violation. 

If $\lfloor q_1 \rfloor =0$ (or similarly for state 2), then $\lceil q_1 \rceil = 1$ and the priority for giving state 1 its first seat is infinite. In that case, the corresponding inequality is vacuously true and the same counting argument above applies. The very large priority ensures state 1 receives its ceiling (=1), and the rest of the argument is unchanged. 

If $q_1$ or $q_2$ are integers, then the proof is similar.
\end{proof}

\begin{remark}[Generalization to $n$ states]
Although Theorem \ref{q.v_criteria_test} is stated for three states, the underlying ideas extend in principle to $n \geq 4$. The main complication already appears at $n=4$: there is no analogue of Theorem \ref{extra} providing a single canonical form for quota violations. Instead, one must account for multiple distinct configurations, each requiring its own criterion, leading to a disjunction of conditions. The number of such cases grows rapidly with $n$, making the resulting criterion unwieldy. For clarity and focus, we restrict attention to the case $n=3$ and leave the general $n$-state setting for future work.
\end{remark}

\subsection{Geometric Visualization of Lower Quota Violations}\label{geo_lower_section}
Using the Lower Quota Violation Criteria Test (\ref{q.v_criteria_test}) and Theorem \ref{l_tri_def}, we can display all quota violations on the simplex of standard quotas $\{(q_1, q_2, q_3) \mid q_1+q_2+q_3=M \}$ within the Adams, Dean, or Huntington-Hill method. This creates a visualization of the geometry and structure of lower quota violations. For instance, below is a simplex plot of all violatory $(q_1, q_2, q_3)$ in the Huntington-Hill method with $M=13$:
\begin{figure}[H]
 \centering
 \includegraphics[width=1\linewidth]{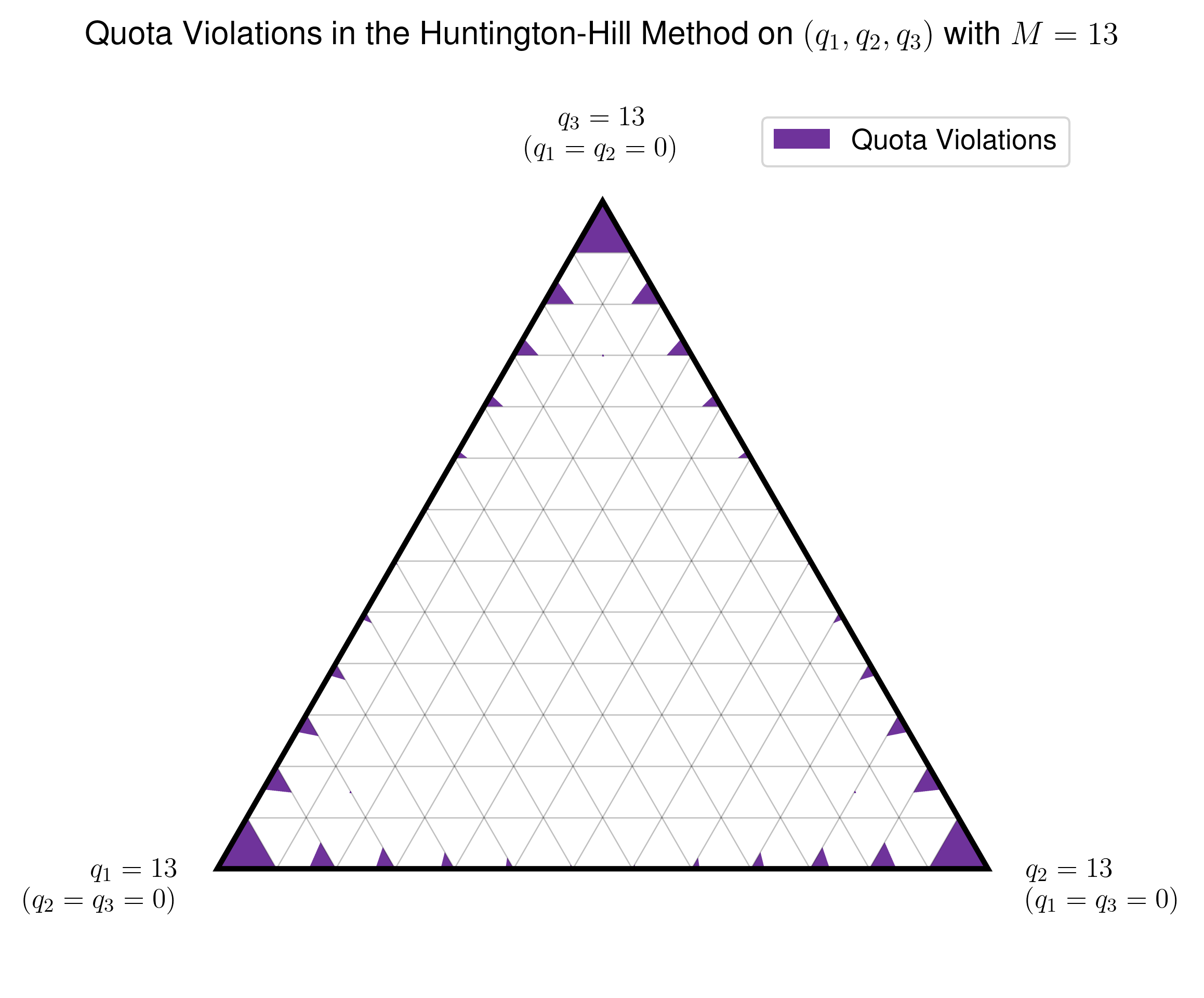}
 \caption{Quota Violations in the Huntington-Hill Method on $\{(q_1, q_2, q_3)\mid q_1+q_2+q_3=M\}$}
 \end{figure}

\subsection{Probability Function for Lower Quota Violations in Three States}
For a fixed $M$, $n=3$ states, and with $A$ the Adams, Dean, or Huntington-Hill method, this section uses the Lower Quota Violation Criteria Test (\ref{q.v_criteria_test}) to find the probability of a lower quota violation. 

Additionally, as the set $\{(q_1,q_2)\,|\, q_1 \text{ or }q_2 \in \mathbb{Z}\}$ has area zero, assume that quotas of the smallest two states, $q_1, q_2$ are not integers. Therefore, Condition (3) of the Lower Quota Violation Criteria Test (\ref{q.v_criteria_test}) is then $\ffo + \fft + \lfloor q_3 \rfloor = M-1$. 

\subsubsection{Motivating Example}
As an example, let $M = 100, \lfloor q_1 \rfloor = 1, \lfloor q_2 \rfloor = 2$, and $A$ the Huntington-Hill method. Write $q_1 = \lfloor q_1 \rfloor +d_1$ and $q_2=\lfloor q_2 \rfloor +d_2$, where $d_1, d_2$ are the decimal parts. The feasible region of the Lower Quota Violation Criterion Test is graphed below:

\begin{figure}[H]
 \centering
 \includegraphics[width=0.8\linewidth]{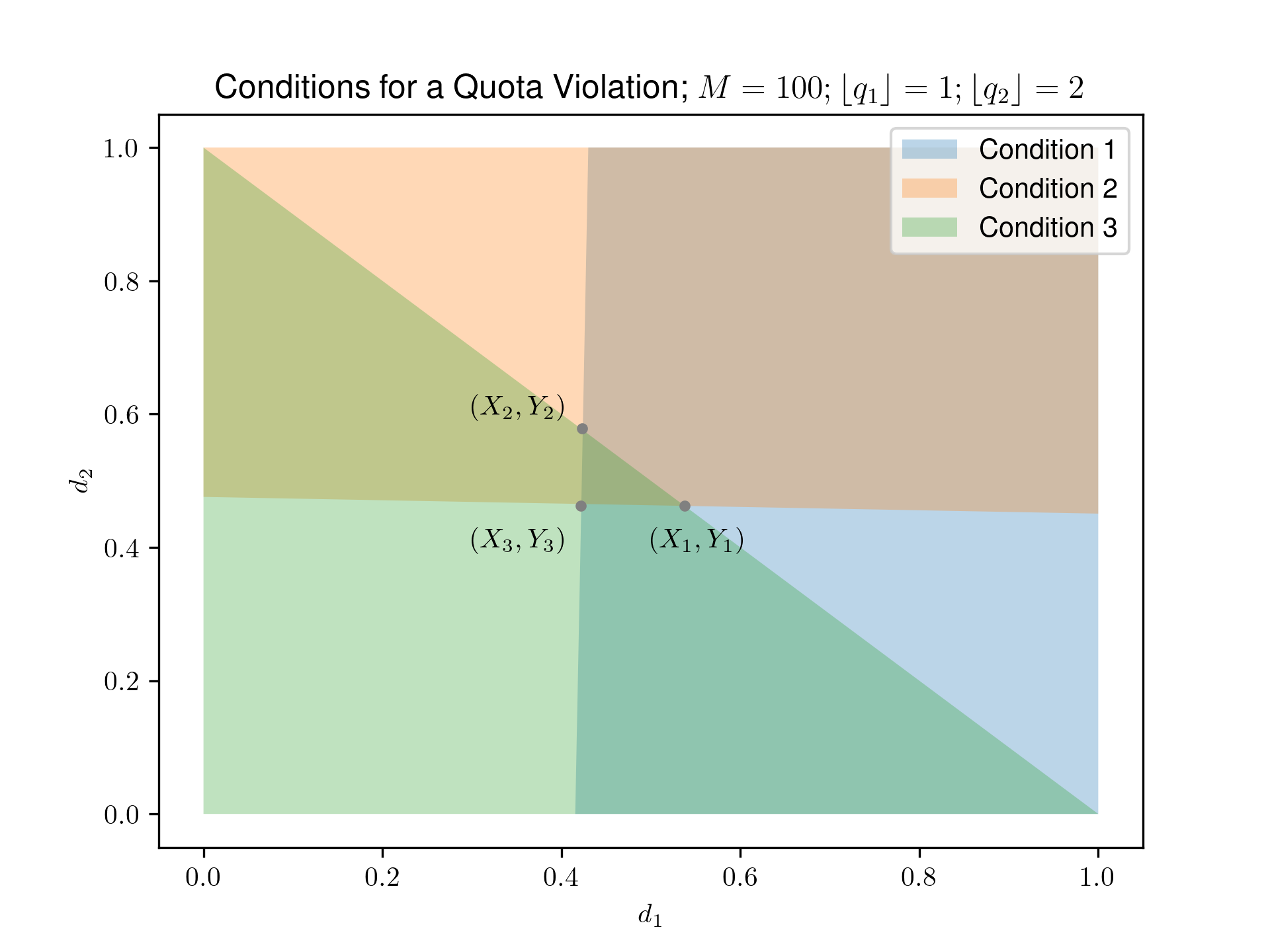}
 \caption{Example Inequality Region with Overlap}
 \end{figure}

The probability of a quota violation in the case that $ \ffo=1$, $\fft=2$, with $q_1$ and $q_2$ assumed to be uniformly distributed, is the area of the feasible region.

Note the feasible region could be empty as in the example with $M = 100, \lfloor q_1 \rfloor = 20, \lfloor q_2 \rfloor = 27$, and $A$ the Huntington-Hill method as shown below: 

\begin{figure}[H]
 \centering
 \includegraphics[width=0.8\linewidth]{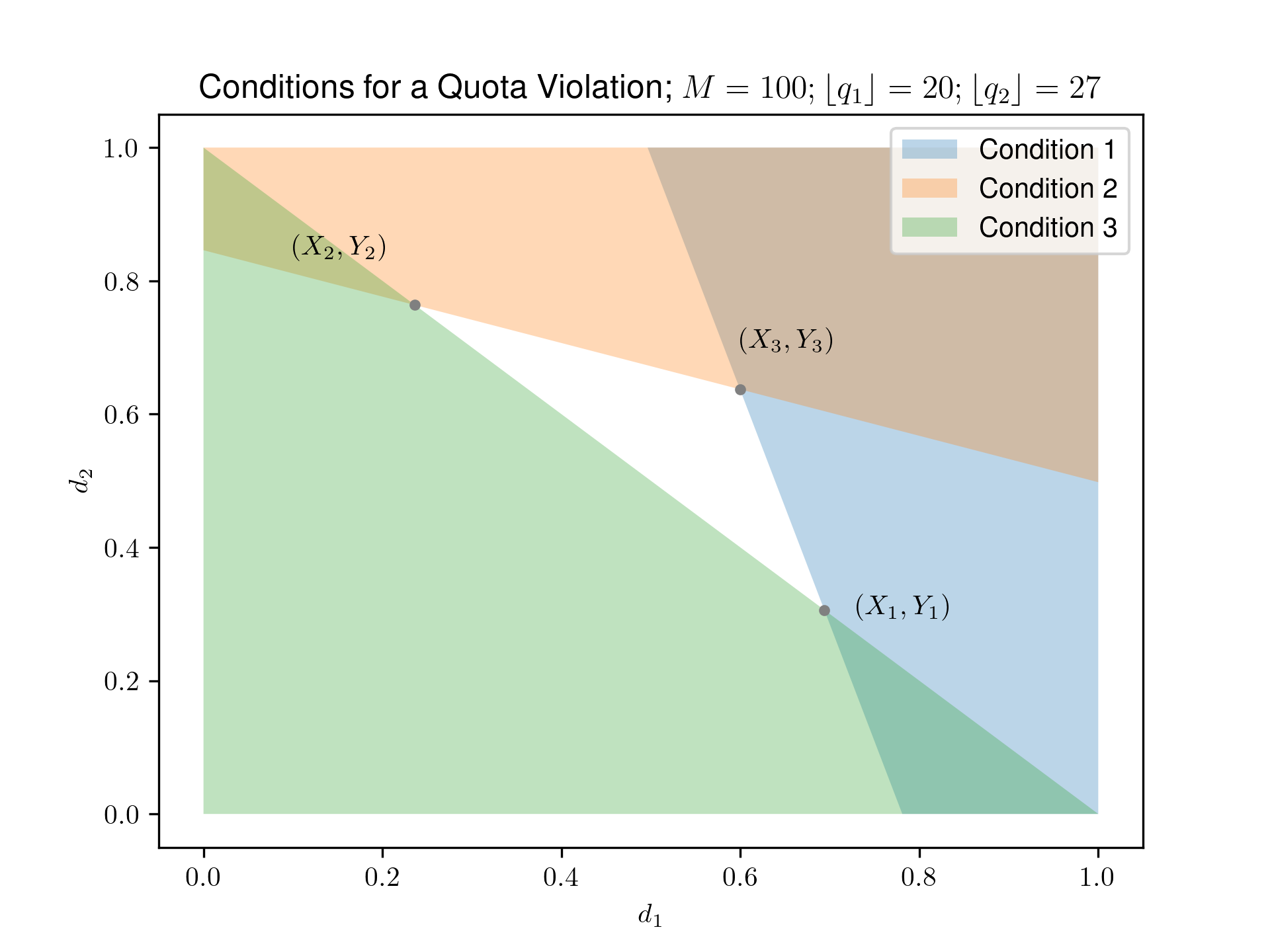}
 \caption{Example Inequality Region without Overlap}
 \end{figure}

\subsubsection{Calculation of the Feasible Region}
The triangular shape of the feasible region in the examples above is not a coincidence. The next theorem shows this is always the case as the feasible region is bounded by three linear inequalities. 

\begin{theorem}\label{l_tri_def} (Equations of the Lower Quota Violation Feasible Region) For $M$ seats, $n = 3$ states, and standard quotas $q_1, q_2$, define $R = M-\lfloor q_1 \rfloor-\lfloor q_2 \rfloor$. Then the feasible region of the Lower Quota Violation Criteria Test (\ref{q.v_criteria_test}) is the subset inside the unit square $(0,1) \times (0,1)$ defined by the following linear inequalities:
$$\begin{array}{clc}
y &>\left(\dfrac{d(R-2)+d(\lfloor q_1 \rfloor)}{-d(\lfloor q_1 \rfloor)}\right)x + \dfrac{Rd(\lfloor q_1 \rfloor)-\lfloor q_1 \rfloor d(R-2)}{d(\lfloor q_1 \rfloor)} & (1)\\ [1.5em]
y &>\left(\dfrac{-d(\lfloor q_2 \rfloor)}{d(R-2)+d(\lfloor q_2 \rfloor )}\right)x +\dfrac{Rd(\lfloor q_2 \rfloor )-\lfloor q_2 \rfloor d(R-2)}{d(R-2)+d(\lfloor q_2 \rfloor)} &(2)\\[1.5em]

y &< 1-x &(3)
\end{array}$$

with division by $0$ interpreted as a vertical line.

\end{theorem}

\begin{proof}
This is a brute force calculation with given floors $\ffo,\fft$. Assume $\ffo, \fft > 0$ and $d(\ffo),d(\fft) \neq 0$. Write the quotas as:$$ (q_1, q_2, q_3) = (\ffo+d_1, \fft+d_2, M-\ffo-\fft-d_1-d_2),$$ where $d_1, d_2 \in [0,1)$ are the decimal parts of $q_1, q_2$. 

Substituting condition (3) of the Lower Quota Violation Criteria Test (\ref{q.v_criteria_test}) gives 
 $$\ffo + \fft + \lfloor M-\ffo-\fft-d_1-d_2 \rfloor = M-1.$$
If $d_1+d_2 \leq 1$, then the floor of the third term is $M-\ffo-\fft-1$ and the equality is satisfied. If $d_1+d_2 > 1 $, the third term is $M-\ffo - \fft -2$ violating condition (3). Thus a lower quota violation is possible only if $d_1+d_2 \leq 1$.

Conditions (1) and (2) of the Lower Quota Violation Criteria Test become:
\[\left(M-\ffo-\fft-d_1-d_2\right)d\left(\ffo\right)<\left(\ffo+d_1\right)d\left(M-\ffo-\fft-2\right),\]
\[\left(M-\ffo-\fft-d_1-d_2\right)d\left(\fft\right)<\left(\fft+d_2\right)d\left(M-\ffo-\fft-2\right).\]

For condition (1), using $R = M-\lfloor q_1 \rfloor-\lfloor q_2 \rfloor$, solve for $d_2$:
\[ d_2>\left(\frac{d(M-R-2)+d(\lfloor q_1 \rfloor)}{-d(\lfloor q_1 \rfloor)}\right)d_1 + \frac{Rd(\lfloor q_1 \rfloor)-\lfloor q_1 \rfloor d(R-2)}{d(\lfloor q_1 \rfloor)}\]

Similarly, solving for $d_2$ in condition (2) gives:
\[d_2>\left(\frac{-d(\lfloor q_2 \rfloor)}{d(R-2)+d(\lfloor q_2 \rfloor )}\right)d_1 +\frac{Rd(\lfloor q_2 \rfloor )-\lfloor q_2 \rfloor d(R-2)}{d(R-2)+d(\lfloor q_2 \rfloor)} \]

The feasible region, when non-empty, is the overlap of the above linear inequalities defined over $(0,1) \times (0,1)$, representing all possible decimal parts of the quotas $q_1,q_2$. 

If $d(\ffo) =0$ or $d(\fft)=0$, then the respective inequality degenerates to a vertical line, but the feasible region remains well defined.
\end{proof}

\begin{remark}\label{implicit} Denote the feasible region defined above as $\bigtriangleup_l(\ffo,\fft)$. It satisfies the following properties:

 \begin{enumerate}
 \item (Subset of the quota feasible region) For the Dean, Adam, or Huntington-Hill methods, and for all pairs ($\ffo$,$\fft$), $\bigtriangleup_l(\ffo,\fft) \subseteq \{ (q_1, q_2) \mid 0 < q_1 < M-q_1-q_2, 0 < q_2 < M-q_1-q_2 \}$. Indeed if a point in $\bigtriangleup_l$ lay outside this set, a lower quota violation would occur for a state that does not have the largest population, contradicting Lemma \ref{vio_on_largest}. 
 
 \item (Triangular shape) For $M>3$ the intersection of the three linear inequalities defining $\bigtriangleup_l(\ffo,\fft)$ is a triangle. This is because the slope of the first line is always less than $-1$, the slope of the second line is always greater than $-1$, and the slope of the third line is always $-1$.
 
 \item (Contained in the Unit Square) If the divisor function $d$ satisfies $\frac{s}{d(s-1)} \geq \frac{s'}{d(s')}$ for all $s, s' \in \mathbb{N}$ and $s' \leq s$, then the region bounded by the linear inequalities is a proper subset of the unit square. Geometrically, this is because the $y$-intercept of the first linear inequality is always greater than $1$ and the $x$-intercept of the second linear inequality is always greater than $1$, ensuring the triangular intersection is inside the unit square whenever it is nonempty. 
 
 Algebraically, the $y$-intercept condition $$\frac{d\left(\lfloor q_1 \rfloor\right)\left(M-\lfloor q_1 \rfloor-\lfloor q_2 \rfloor\right)-\lfloor q_1 \rfloor d\left(M-\lfloor q_1 \rfloor-\lfloor q_2 \rfloor-2\right)}{d\left(\lfloor q_1 \rfloor\right)} <1$$ rearranges to $$\frac{\ffo}{d(\ffo)} > \frac{\lfloor q_3 \rfloor}{d(\lfloor q_3 \rfloor -1)},$$ where then $M - \ffo - \fft = \lfloor q_3 \rfloor+1$. By symmetry, the same holds for the $x$-intercept. 
 In particular, this holds in Adam's method, Dean method, and the Huntington-Hill Method.
\end{enumerate}
These properties are used implicitly in Corollary \ref{l_uniform_calc}.
\end{remark}

\subsubsection{Probability Formula for Lower Quota Violations}
The probability of a lower quota violation of $A(q_1,q_2,q_3)$ for a fixed $q_1, q_2$ is then the double integral over the feasible region of the joint pdf. Summing this probability for for all possible values of $\ffo,\fft$ gives the following formula: 

\begin{theorem}\label{lower_prob} (Probability of Lower Quota Violation) \label{Probability_lq_updated}
 For $M$ seats and $n=3$ states, let $q_1$ and $q_2$ be the standard quotas of states 1 and 2 with joint pdf $f$, and define $q_3=M-q_1-q_2$. For $q_1, q_2 < q_3$ and for apportionment method $A$ equal to Adams, Dean, or Huntington-Hill, the probability of a lower quota violation is given by:
\[ \sum_{\ffo,\fft} \iint_{\bigtriangleup_l(\ffo,\fft)} f(\ffo+x,\fft+y)\,dx\,dy \]
where $\bigtriangleup_l(\ffo,\fft)$ denotes the feasible region of the Lower Quota Violation Criteria Test (\ref{q.v_criteria_test}) defined in Theorem \ref{l_tri_def}.
\end{theorem}

\begin{corollary}\label{l_uniform_calc}
Let $(q_1, q_2)$ be chosen uniformly from $\{ (q_1, q_2) \mid 0 < q_1 < M-q_1-q_2, 0 < q_2 < M-q_1-q_2 \}$. Then, for the Dean, Adams, and Huntington-Hill method with $M>3$, the probability of a quota violation is:

$$ \frac{6}{M^2}\sum_{\ffo, \fft} \max_{} \LL\{ 0, \LL( \frac{1}{2} \det \begin{bmatrix} \max_{i=1,2}X_i(\ffo,\fft) & 1 - \max_{i=1,2}X_i(\ffo,\fft) & 1 \\ \min_{i=1,2}X_i(\ffo,\fft) & 1 - \min_{i=1,2}X_i(\ffo,\fft) & 1 \\ X_3(\ffo,\fft) & Y_3(\ffo,\fft) & 1 \end{bmatrix} \RR) \RR\},$$

Where, for $R =M-{\ffo}-{\fft}$, $X_i$ and $Y_i$ are defined by:

$$\begin{array}{cl}
 X_1  & = \dfrac{d\left({\fft}\right)+d\left(R-2\right)-Rd\left({\fft}\right)+{\fft}d\left(R-2\right)}{d\left(R-2\right)} \\[1em]
 Y_1  & = \displaystyle{1 - X_1} \\[1em]
 Y_2  & = \dfrac{d({\ffo})+d(R-2)-Rd({\ffo})+{\ffo}d(R-2)}{d(R)}\\
 X_2  & = \displaystyle{1 - Y_2}\\[1em]
 X_3  & = \dfrac{{\fft}d({\ffo})+Rd({\ffo})-{\ffo}d({\fft})-{\ffo}d(R-2)}{d(R-2)+d({\fft})+d({\ffo})}\\[1em]
 Y_3  &= \dfrac{{\ffo}d({\fft})+Rd({\fft})-\fft d({\ffo})-{\fft}d(R-2)}{d(R-2)+d({\fft})+d\left({\ffo}\right)} .
\end{array}$$

Additionally, if $M=3$, the probability of a lower quota violation is $\frac{1}{3}$.
\end{corollary}

\begin{proof}
For each fixed $\ffo,\fft$, the area of $\bigtriangleup_l(\ffo,\fft)$ under a uniform measure is $\iint_{\bigtriangleup_l(\ffo,\fft)} f(\ffo+x, \fft+x) \,dx\,dy$ . The vertices $(X_i, Y_i)$ of $\bigtriangleup_l(\ffo,\fft)$ are given by the formulas above. The area of $\bigtriangleup_l(\ffo,\fft)$ is then:

$$ \max_{} \LL\{ 0, \LL( \frac{1}{2} \det \begin{bmatrix} \max_{i=1,2}X_i(\ffo,\fft) & 1 - \max_{i=1,2}X_i(\ffo,\fft) & 1 \\ \min_{i=1,2}X_i(\ffo,\fft) & 1 - \min_{i=1,2}X_i(\ffo,\fft) & 1 \\ X_3(\ffo,\fft) &Y_3 (\ffo,\fft) & 1 \end{bmatrix} \RR) \RR\}. $$

Where the right hand rule and the location of $(X_3, Y_3)$ determines whether the region is nonempty. Since the total area of $\{ (q_1, q_2) \mid q_1, q_2< M-q_1 -q_2 \}$ is $\frac{M^2}{6}$, the stated probability formula follows.

For $M=3$, in all three methods 
$$\mathrm{Area}( \bigtriangleup_l(0,1))=\mathrm{Area}( \bigtriangleup_l(1,0)) = 0, \mathrm{Area}( \bigtriangleup_l(0,0)) = \frac{1}{2}.$$ The total probability is then $\frac{1}{2} \cdot \frac{6}{3^2} = \frac{1}{3}$.
\end{proof}

A computer program to evaluate the above formula is provided at \url{https://github.com/TylerCWunder/Probability_and_Structure_of_Quota_Violations_in_Divisor_Methods_of_Apportionment_code.git}. 

\subsubsection{Empirical Results}\label{sample1section}
Given an apportionment method $A$, $M$ seats, and $n=3$ states, the probability of a lower quota violation under a uniform distribution is computed using Corollary \ref{l_uniform_calc}. The resulting probabilities are shown in the tables below and compared with a simulation of 100,000 random $(q_1, q_2)$ from a uniform distribution on $\{ (q_1, q_2) \mid 0 < q_1 < M-q_1-q_2, 0 < q_2 < M-q_1-q_2 \}$:
\setlength{\tabcolsep}{6pt}
\renewcommand{\arraystretch}{1.33} 
\begin{center}
\begin{tabular}{|c|c|c|c|}
 \hline
 \multicolumn{4}{|c|}{\textbf{Quota Violations in the Huntington-Hill Method with Uniform Quotas}}\\
 \hline
 
 $M$ & Theoretical Probability & Sample Probability & $95\%$ Confidence Interval \\ \hline
  $3$ & $0.33333$ & $0.33533$ & $\left(0.33240,0.33826\right)$\\
 $5$&$0.13092$ &$0.13055$ &$\left(0.12846,0.13264\right)$ \\
 $10$&$0.04904$&$0.04915$&$\left(0.04781,0.05049\right)$ \\
 $15$&$0.02916$&$0.02906$&$\LL( 0.02802 , 0.03010 \RR)$ \\ 
 $20$ & $0.02047$ & $0.02069$ & $\left(0.01981,0.02157\right)$ \\
 $100$ & $0.00328$ & $0.00351$ &$\left(0.00314,0.00388\right)$ \\
 \hline
\end{tabular}
\end{center}
\begin{center}
\begin{tabular}{|c|c|c|c|}
 \hline
 \multicolumn{4}{|c|}{\textbf{Quota Violations in the Adams Method with Uniform Quotas}}\\
 \hline
$M$ & Theoretical Probability & Sample Probability & $95\%$ Confidence Interval \\ \hline
  $3$ & $0.33333$ & $0.33249$ & $( 0.32957 , 0.33541 )$\\
 $5$&$0.16000$ &$0.15852$ &$( 0.15626 , 0.16078 )$ \\
 $10$&$0.11057$&$0.1099$&$( 0.10796 , 0.11184 )$ \\
 $15$&$0.09897$&$0.09961$&$( 0.09775 , 0.10147 )$ \\ 
 $20$ & $0.09365$ & $0.09143$ & $( 0.08964 , 0.09322 )$ \\
 $100$ & $0.08210$ & $0.08188$ &$( 0.08018 , 0.08358 )$ \\
\hline
\end{tabular}
\end{center}
\begin{center}
\begin{tabular}{|c|c|c|c|}
 \hline
 \multicolumn{4}{|c|}{\textbf{Quota Violations in the Dean Method with Uniform Quotas}}\\
 \hline
$M$ & Theoretical Probability & Sample Probability & $95\%$ Confidence Interval \\ \hline

  $3$ & $0.33333$ & $0.33363$ & $( 0.33071 , 0.33655 )$\\
 $5$&$0.13714$ &$0.13818$ &$( 0.13604 , 0.14032 )$ \\
 $10$&$0.05581$&$0.05545$&$( 0.05403 , 0.05687 )$ \\
 $15$&$0.03397$&$0.03412$&$( 0.03299 , 0.03525 )$ \\ 
 $20$ & $0.02403$ & $0.02444$ & $( 0.02348 , 0.02540 )$ \\
 $100$ & $0.00376$ & $0.00362$ &$( 0.00325 , 0.00399 )$ \\
\hline
\end{tabular}
\end{center}
The code used for random sampling is available at \url{https://github.com/TylerCWunder/Probability_and_Structure_of_Quota_Violations_in_Divisor_Methods_of_Apportionment_code.git}.

Comparisons to simulations under different distributions may be found in \hyperref[AppA]{Appendix A}.

\section{Structure and Probability of Upper Quota Violations in the Jefferson Method for Three States}
\subsection{Analogous Theoretical Results}
We now provide a similar study of upper quota violations in the Jefferson method when apportioning three states. The Webster, Adams, Dean, and the Huntington-Hill methods have no upper quota violations by Lemma \ref{no-up} so only the Jefferson method needs to be considered. 

\begin{lemma}\label{up_3_structure}
In a divisor apportionment method with $M$ seats, $n=3$ states, and standard quotas $q_i$, if an upper quota violation occurs then the apportionment of the state with the violation is $\lceil q_i \rceil +1$, and the other states are assigned their lower quotas.
\end{lemma}\begin{proof}
Assume without loss of generality that the upper quota violation occurs on state 3. The claim is that the final apportionment is then $(a_1,a_2,a_3)=(\ffo, \fft, \lceil q_3 \rceil +1)$. Since the first two states cannot be apportioned less than the floor of the quota by Theorem \ref{no_up_and_low}, $a_1 \geq \ffo$ and $a_2 \geq \fft$. Additionally, as an upper quota violations is assume $a_2 \geq \lceil q_3 \rceil +1$.

Write $q_i = \lfloor q_i \rfloor +d_i$, where $d_i \in [0,1)$ is the decimal part. Then, in particular $\sum d_i \in \{0,1,2\}$.

Then, assume $d_3 \neq 0$. As such, $\lfloor q_3 \rfloor = \lceil q_3 \rceil -1$ and
\[ a_1 + a_2 = M - a_3 \leq M - \lceil q_i \rceil -1 =\sum \lfloor q_i \rfloor + \sum d_i - \lceil q_i \rceil -1 \leq \ffo + \fft + \sum d_i -2 . \]
Then, for each case of $\sum d_i \in \{0,1,2\}$:

If $\sum d_i = 0$, then all of the quotas are integers and the apportionment is not a quota violation.

If $\sum d_i = 1$, then $a_1 + a_2 \leq \ffo + \fft - 1$, implying that one of states one or two is assigned less than their lower quota, which is a contradiction with the earlier statement $a_1 \geq \ffo$ and $a_2 \geq \fft$.

If $\sum d_i = 2$, then $a_1 + a_2 \leq \ffo + \fft$. As it is also the case that $a_1 \geq \ffo$ and $a_2 \geq \fft$, this suggests that they are all in fact equal. Then, $a_3 = M - a_1 - a_2 = \sum \lfloor q_i \rfloor + \sum d_i - \ffo - \fft = \lfloor q_3 \rfloor + 2 = \lceil q_3 \rceil +1$. Thus, $(a_1,a_2,a_3)=(\ffo, \fft, \lceil q_3 \rceil +1)$ and as this is the only possible case for $d_3 \neq 0$ the claim is proven under this assumption.

Then, if $d_3 = 0$:
\[ a_1 + a_2 = M - a_3 \leq M - \lceil q_i \rceil -1 =\sum \lfloor q_i \rfloor + \sum d_i - \lceil q_i \rceil -1 \leq \ffo + \fft + \sum d_i -1 . \]
Which gives for each case of $\sum d_i \in \{0,1,2\}$:

If $\sum d_i = 0$, then all of the quotas are integers and the apportionment is not a quota violation.

If $\sum d_i = 2$, then as $d_i \in [0,1)$ this violates the assumption that $d_3=0 $.

If $\sum d_i = 1$, then $a_1 + a_2 \leq \ffo + \fft$, which grants $a_1 = \lfloor q_1 \rfloor$ and $a_2 = \lfloor q_2 \rfloor$. Then, $a_3 = M - a_1 - a_2 = \sum \lfloor q_i \rfloor + \sum d_i - \ffo - \fft = \lfloor q_3 \rfloor + 1 = \lceil q_3 \rceil +1$. Thus, $(a_1,a_2,a_3)=(\ffo, \fft, \lceil q_3 \rceil +1)$ and as this is the only possible case for $d_3 = 0$ the claim is proven here as well.
\end{proof}

Next, similar to Lemma \ref{vio_on_largest}, we justify using $q_3$ in the above proof as having the upper quota violation. 

\begin{lemma}\label{up_vioonlarge}
In the Jefferson method and $n=3$ states, upper quota violations can only occur on the state with the largest population.
\end{lemma}

\begin{proof}
Order the state's populations $p_1 < p_2 < p_3$. Assume there is an upper quota violation on state $i$ with $i \in \{1,2\}$. Lemma \ref{up_3_structure} then implies both that state three is not apportioned $\lceil q_3 \rceil$ as state $i$ is apportioned $\lceil q_i \rceil +1.$ This gives the following relation between priority values:
\[ \frac{q_i}{d( \lceil q_i \rceil)} >\frac{q_3}{d( \lfloor q_3 \rfloor)} \]
Rearranging gives:
$$ q_i ( \lfloor q_3 \rfloor+1) > q_3 ( \lceil q_i \rceil + 1)$$
The next Lemma will show that this implies $q_i > q_3$, a contradiction.
\end{proof}

\begin{lemma} \label{ineq2}
For $a,b \in (0, \infty)$, $a ( \lfloor b \rfloor+1) > b ( \lceil a \rceil + 1)$ implies $a>b$. \end{lemma}
\begin{proof}
If $b \in \mathbb{Z}$, then the statement is $a (b +1) > b( \lceil a \rceil + 1)$, which rearranges to $ab + a >b \lceil a \rceil +b \geq ab + b$, hence $a > b$. If $b$ is not an integer, the statement rearranges to $a \lceil b \rceil >b \lceil a \rceil +b \geq \lfloor b \rfloor a + b$. Then, $b<a ( \lceil b \rceil - \lfloor b \rfloor) = a$.
\end{proof}

Summarizing Lemma \ref{up_3_structure} and \ref{up_vioonlarge} gives the following:
\begin{theorem}\label{extra2} (Classification of Upper Quota Violations in 3 States) 
Order the populations $p_1 < p_2<p_3$. If $A$ is the Jefferson method and $A(p_1, p_2, p_3)$ has a quota violation, the final apportionment must be $(a_1,a_2,a_3)= (\lfloor q_1 \rfloor, \lfloor q_2 \rfloor, \lceil q_3 \rceil +1)$, where $q_i$ is the standard quota of state $i$. That is, the quota violation must be an upper quota violation occurring in the third state. 
\end{theorem}

\subsubsection{Upper Quota Violation Criteria Test for 3 States}
Next a test is provided to check for an upper quota violation on three states under Jefferson Apportionment. Similar to Theorem \ref{q.v_criteria_test}, this provides an analytic and more direct way to check for a quota violation without needing to run the apportionment algorithm. This test will then be used to find the probability of an upper quota violation for 3 states.

\begin{theorem} \label{uq.v_criteria_test} \textit{(Upper Quota Violation Criteria Test)}
Let $A$ be the Jefferson method, with divisor function $d(s)$. Order the populations such that $p_1 <p_2 < p_3$ and calculate the standard quotas $q_i$. Then the apportionment $A(p_1, p_2, p_3)$ on $M$ seats has an upper quota violation if and only if all three of the following statements are true:
\begin{enumerate}
\item $q_3 d( \lfloor q_1 \rfloor) > q_1 d(\lceil q_3 \rceil)$
\item $q_3 d( \lfloor q_2 \rfloor) > q_2 d(\lceil q_3 \rceil)$
\item $\lfloor q_1 \rfloor + \lfloor q_2 \rfloor + \lceil q_3 \rceil = M-1$ (equivalently for $q_3 \notin \mathbb{Z}$, $\ffo + \fft + \lfloor q_3 \rfloor = M-2$.)

\end{enumerate}
\end{theorem}
The proof is similar to the proof of Theorem \ref{q.v_criteria_test}.
\begin{proof}
Assume that there is an upper quota violation. Then, by Theorem \ref{extra2}, the apportionment is $(a_1, a_2, a_3) = (\lfloor q_1 \rfloor, \lfloor q_2 \rfloor, \lceil q_3 \rceil +1)$. Then, states one and two are now assigned the floors of their quotas and state three is assigned the ceiling of its quota plus one. Writing the relationship between priority values suggested by this grants conditions (1) and (2). For condition (3), $M = a_1 + a_2 + a_3 = \ffo + \fft +\lceil q_3 \rceil +1$.

Assume that the three conditions hold. Then, calculate the apportionment by first assigning each state $\lfloor q_i \rfloor$, where every state is guaranteed at least that many seats by Theorem \ref{nvv}. Then, assign the remaining seats by following the priority value relations implied by conditions (1) and (2), which is assigning all additional seats to state three. This creates the apportionment $(a_1, a_2, a_3) = (\lfloor q_1 \rfloor, \lfloor q_2 \rfloor, \lceil q_3 \rceil +1)$, which is an upper quota violation.
\end{proof}

From there, in analogy with Theorem \ref{l_tri_def}, derive the following result:
\begin{theorem}\label{u_tri_def} (Equations of the Upper Quota Violation Feasible Region) For $M$ seats, $n = 3$ states, and standard quotas $q_1, q_2$, define $R = M-\lfloor q_1 \rfloor-\lfloor q_2 \rfloor$. Then the feasible region of the Upper Quota Violation Criteria Test (\ref{uq.v_criteria_test}) is the subset inside the unit square $(0,1) \times (0,1)$ defined by the following linear inequalities:
$$\begin{array}{clc}
\displaystyle{y} &< \displaystyle{\LL(\dfrac{d(R-1)+d(\lfloor q_1 \rfloor)}{-d(\lfloor q_1 \rfloor)} \RR)x +\space\dfrac{Rd(\lfloor q_1 \rfloor)-\lfloor q_1 \rfloor d(R-1)}{d\left(\lfloor q_1 \rfloor\right)}} &(1)\\[1.5em]
y &< \displaystyle{ \left(\frac{-d(\lfloor q_2 \rfloor)}{d(R-1)+d(\lfloor q_2 \rfloor)}\right)x + \frac{Rd(\lfloor q_2 \rfloor)-\lfloor q_2 \rfloor d(R-1)}{d(R-1)+d(\lfloor q_2 \rfloor)}} &(2)\\[1.5em]
y &> \displaystyle{1-x} &(3)
\end{array}$$

where division by $0$ is interpreted as a vertical line.
\end{theorem}

The proof is omitted as it is similar to the proof of Theorem \ref{l_tri_def}. Similar properties to Remark \ref{implicit} apply to this region.
\subsubsection{Geometric Visualization of Upper Quota Violations}
Similarly to Section \ref{geo_lower_section}, one may use Theorem \ref{u_tri_def} to display all violatory $(q_1, q_2, q_3)$ Un the Jefferson Method on the simplex $\{(q_1, q_2, q_3) \mid q_1+q_2+q_3=M \}$. This allows for a visualization of the geometry and structure of upper quota violations. For instance, below is a simplex plot indicating all violatory $(q_1, q_2,q_3)$ in the Jefferson method with $M=13$:

\begin{figure}[H]
 \centering
 \includegraphics[width=1\linewidth]{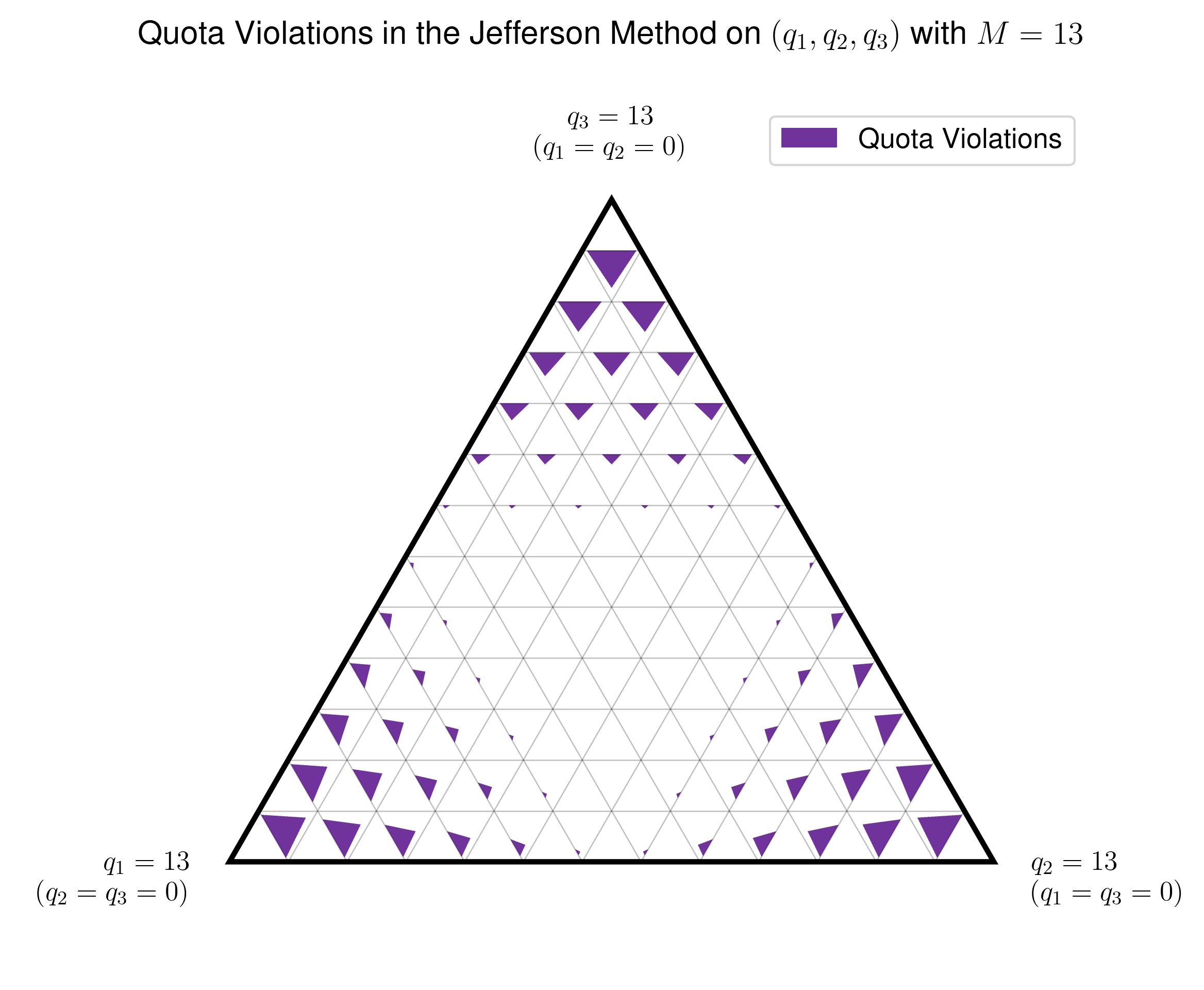}
 \caption{Quota Violations in the Jefferson Method on $\{(q_1, q_2, q_3)\mid q_1+q_2+q_3=M\}$}
 \end{figure}

\subsubsection{Probability Formula for Upper Quota Violations}
For a fixed number of seats $M$, the probability of an upper quota violation of $A(q_1, q_2, q_3)$ for a fixed $q_1, q_2$ is then the double integral over the feasible region of the joint pdf. Summing this probability for for all possible values of $\ffo, \fft$ gives the following formula:

\begin{theorem}\label{up_prob} (Probability of Upper Quota Violation) \label{Probability_updated}
 For $M$ seats and $n=3$ states, let $q_1$ and $q_2$ be the standard quotas of states 1 and 2 with joint pdf $f$, and define $q_3=M-q_1-q_2$. For $q_1, q_2 < q_3$ and the Jefferson method, the probability of an upper quota violation is given by:
 
\[ \sum_{\ffo,\fft} \iint_{\bigtriangleup_u(\ffo,\fft)} f(\ffo+x,\fft+y)\,dx\,dy \]
where $\bigtriangleup_u(\ffo,\fft)$ denotes the feasible region of the Upper Quota Violation Criteria Test (\ref{uq.v_criteria_test}) defined in Theorem \ref{u_tri_def}.
\end{theorem}
\begin{corollary}\label{up_uni_prob}
If $(q_1, q_2)$ are chosen uniformly on $\{ (q_1, q_2) \mid 0 < q_1 < M-q_1-q_2, 0 < q_2 < M-q_1-q_2 \}$, then for the Jefferson method the probability of an quota violation for any $M$ is given by: 
\[ \frac{6}{M^2}\sum_{\ffo, \fft} \max_{} \LL\{ 0, \LL( \frac{1}{2} \det \begin{bmatrix} \min_{i=1,2}X_i(\ffo,\fft) & 1 - \min_{i=1,2}X_i(\ffo,\fft) & 1 \\ \max_{i=1,2}X_i(\ffo,\fft) & 1 - \max_{i=1,2}X_i(\ffo,\fft) & 1 \\ X_3(\ffo,\fft) & Y_3(\ffo,\fft) & 1 \end{bmatrix} \RR) \RR\}. \]
Where $X_i$ and $Y_i$ are defined, with $R=M-\ffo-\fft$, as:
$$\begin{array}{cl}
X_1 &= \dfrac{Rd({\ffo})-{\ffo}d(R-1)-d({\ffo})}{d(R-1)} \\[1em]
Y_1 &= \displaystyle{1 - X_1}\\[1em]
Y_2 &= \dfrac{Rd({\fft})-{\fft}d(R-1)-d({\fft})}{d(R-1)} \\[1em]
X_2 &= \displaystyle{1 - Y_2}\\[1em]
X_3 &= \dfrac{Rd({\ffo})-{\ffo}d(R-1)-{\ffo}d({\fft})+{\fft}d({\ffo})}{d(R-1)+d({\ffo})+d({\fft})}\\[1em]
Y_3 &= \dfrac{Rd({\fft})-{\fft}d(R-1)-{\fft}d({\ffo})+{\ffo}d({\fft})}{d(R-1)+d({\ffo})+d({\fft})}.
\end{array}$$
\end{corollary}
The proof is similar to that of Corollary \ref{l_uniform_calc}.

A computer program to evaluate the above formula is provided at \url{https://github.com/TylerCWunder/Probability_and_Structure_of_Quota_Violations_in_Divisor_Methods_of_Apportionment_code.git}.

\subsubsection{Empirical Results}\label{sample2section}
Given the Jefferson method with $M$ seats and $n=3$ states, the probability of an upper quota violation under a uniform distribution is computed using Corollary \ref{up_uni_prob}. The resulting probabilities are shown in the tables below and compared with a simulation of 100,000 random $(q_1, q_2)$ from a uniform distribution on $\{ (q_1, q_2) \mid 0 < q_1 < M-q_1-q_2, 0 < q_2 < M-q_1-q_2 \}$:
\begin{center}
\begin{tabular}{|c|c|c|c|}
 \hline
 \multicolumn{4}{|c|}{\textbf{Quota Violations in the Jefferson Method with Uniform Quotas}}\\
 \hline
$M$ & Theoretical Probability & Sample Probability & $95\%$ Confidence Interval \\ \hline
 $3$&$0.02222$ & $0.02279$ & $( 0.02187 , 0.02371 )$\\
 $5$&$0.03943$ &$0.03888$ &$( 0.03768 , 0.04008 $ \\
 $10$&$0.05671$&$0.05694$&$( 0.05550 , 0.05838 )$ \\
 $15$&$0.06362$&$0.06114$&$( 0.05966 , 0.06262 )$ \\ 
 $20$ & $0.06729$ & $0.06693$ & $( 0.06538, 0.06848)$ \\
 $100$ & $0.07686$ & $0.07548$ &$( 0.07384 , 0.07712 )$ \\
 
 \hline
\end{tabular}
\end{center}
The code used for random sampling is available at \url{https://github.com/TylerCWunder/Probability_and_Structure_of_Quota_Violations_in_Divisor_Methods_of_Apportionment_code.git}.

Comparisons to simulations under different distributions may be found in \hyperref[AppA]{Appendix A}.
\section{Discussion of Empirical Results}
The empirical results presented in Sections \ref{sample1section}, \ref{sample2section}, and \hyperref[AppA]{Appendix A} provide strong validation of the theoretical probability formulas derived in Theorems \ref{Probability_lq_updated} and \ref{Probability_updated}. Across all tested apportionment methods (Adams, Dean, Huntington-Hill, and Jefferson), and across varying values of $M$, the theoretical probabilities fall consistently within the $95\%$ confidence intervals of the simulated sample probabilities, with 100,000 random draws per configuration. This agreement holds not only under the uniform distribution on the simplex (equivalently, Dirichlet $(1,1,1)$ and IID exponential populations) but also under the Dirichlet $(2,2,2)$ and Dirichlet $(\frac12 , \frac12 , \frac12)$ distributions, as well as uniform populations on an interval. 

Several patterns emerge from the data: For Adams method, violation probabilities decrease only modestly as $M$ grows, falling from roughly $33\%$ at $M = 3$ to about $8\%$ at $M = 100$ under the uniform distribution. This reflects the method's persistent tendency to round up aggressively and thereby favor smaller states. The Huntington-Hill and Dean methods exhibit more rapid decay in violation probability with increasing M, suggesting their rounding rules become better aligned with proportionality as the number of seats grows. The Jefferson method, by contrast, shows the opposite trend: upper quota violation probabilities increase with $M$, approaching roughly $8\%$ at $M = 100$, consistent with Jefferson's well-known bias toward larger states. 

Additionally, the decrease in quota violations from the Dirichlet $(\frac12 , \frac12 , \frac12)$ to the Dirichlet $(1,1,1)$ and Dirichlet $(2,2,2)$ distributions on the Huntington-Hill, Adams, and Dean methods for all values of $M$ tested and in the Jefferson method with $M=15$ highlights how quota violations are more prevalent towards the edge of the simplex.

These patterns across methods reinforce the theoretical classification results and further highlight the underlying geometry of quota violations.

\section{Conclusion and Next Steps}

This work provides a comprehensive analysis of quota violations in divisor methods for three states, advancing the quantitative and qualitative understanding of the Balinski–Young Theorem. The main contribution of the paper lies in making the consequences of the Balinski-Young theorem explicit and quantitative in this restricted setting by providing a complete classification of possible quota violations and analyzing the conditions and frequencies under which they occur. The work thus refines and illustrates a classical impossibility result rather than challenging it or proposing an alternative apportionment method

The results are significant because they illustrate how probabilistic and geometric techniques are combined to analyze apportionment. Specifically, these methods both quantify the frequency and illuminate the structure of quota violations in divisor methods. The ability to move beyond impossibility theorems toward concrete probability calculations represents a step forward in the applied probability analysis of voting and apportionment rules.

Several natural directions for future research emerge from this study. First and foremost would be the generalization to a larger numbers of states and all divisor methods, where the geometry of quota violations becomes richer and more intricate. Second, extending the framework beyond divisor methods to rank-index methods or even to the entire class of apportionment methods covered by the Balinski-Young Theorem would provide a unified probabilistic perspective on fairness violations in apportionment. Finally, a computational analysis where the number of seats, states or ratios between populations grows large could offer insights into the limiting behavior of quota violation probabilities.

The results presented here make a meaningful contribution to the theory and practice of apportionment. They not only clarify the mathematical structure of quota violations in divisor methods but also open new avenues for probabilistic, geometric, and computational investigations. We believe that the framework developed here will serve as a foundation for further research at the intersection of applied probability, social choice theory, and political representation.

\section{Declarations}
\subsection{Ethical Approval and Consent to Participate}
Not applicable. 

\subsection{Competing Interests} The authors declare that they have no known competing financial interests or personal relationships that could have appeared to influence the work reported in this paper.

\subsection{Data Availability Statement} The raw data generated during this study is available upon request. All relevant processed data is included within this article and the appendix. The Python code is available to run
here: \url{https://github.com/TylerCWunder/Probability_and_Structure_of_Quota_Violations_in_Divisor_Methods_of_Apportionment_code.git}

\subsection{Consent for publication}
Not applicable.

\subsection{Funding}
This research was supported by Johns Hopkins’
Provost’s Undergraduate Research Award (PURA)

\newpage
\nocite{huntington-a-new}
\nocite{ElHelaly2019}
\printbibliography

\begin{contact}
Joseph W. Cutrone\\
Department of Mathematics\\
Johns Hopkins University\\
Baltimore, Maryland, USA\\
\email{jcutron2@jhu.edu}\\\\

Tyler C. Wunder\\
Johns Hopkins University\\
Baltimore, Maryland, USA\\
\email{twunder2@jh.edu}
\end{contact}

\appendix
\newpage

\section{Comparison of Theoretical Results to Simulations Under Various Distributions}\label{AppA}
The following are the results of additional simulations and comparisons to theoretical values calculated using Theorems \ref{lower_prob} and \ref{up_prob}.

The python code used to sample quota violations is available here: \url{https://github.com/TylerCWunder/Probability_and_Structure_of_Quota_Violations_in_Divisor_Methods_of_Apportionment_code.git}

\begin{center}
\begin{tabular}{|c|c|c|c|c|}
 \hline
 \multicolumn{5}{|c|}{\textbf{Probability of Quota Violations with $(\frac{q_1}{M}, \frac{q_2}{M}, \frac{q_3}{M}) \sim \mathrm{Dir}(\frac{1}{2},\frac{1}{2},\frac{1}{2})$}}\\
 \hline
 
 Method & $M$ & Theoretical Probability & Sample Probability & $95\%$ Confidence Interval \\ \hline
 Huntington-Hill & $3$ & $0.55051$ & $0.5527$ & $(0.54962,0.55578)$ \\
 Adams & $3$ & $0.55051$ & $0.55421$ & $(0.55113,0.55729)$ \\
 Dean & $3$ & $0.55051$ & $0.54921$ & $(0.54613,0.55229)$ \\
 Jefferson & $3$ & $0.01252$ & $0.01283$ & $(0.01213,0.01353)$ \\
 Huntington-Hill & $5$ & $0.33758$ & $0.33468$ & $(0.33176,0.33760)$ \\
 Adams & $5$ & $ 0.37385$ & $0.37233$ & $ ( 0.36933 , 0.37533 )$ \\
 Dean & $5$ & $0.34624$ & $0.34727$ & $(0.34432,0.35022)$ \\
 Jefferson & $5$ & $0.03003$ & $0.02934$ & $(0.02829,0.03039)$ \\
 Huntington-Hill & $10$ & $0.20313$ & $0.20478$ & $(0.20228,0.20728)$ \\
 Adams & $10$ & $0.29367$ & $0.29473$ & $(0.29190,0.29756)$ \\
 Dean & $10$ & $0.21604$ & $0.21577$ & $(0.21322,0.21831)$ \\
 Jefferson & $10$ & $0.05637$ & $0.05682$ & $(0.05539,0.05825)$ \\
 Huntington-Hill & $15$ & $0.15545$ & $0.15436$ & $( 0.15212 , 0.15660 )$ \\
 Adams & $15$ & $0.26701$ & $0.26744$ & $(0.26470,0.27018)$ \\
 Dean & $15$ & $0.16703$ & $0.1651$ & $(0.16280,0.16740)$ \\
 Jefferson & $15$ & $0.07108$ & $0.07055$ & $(0.06896,0.07214)$ \\

 \hline
\end{tabular}
\end{center}

\begin{center}
\begin{tabular}{|c|c|c|c|c|}
 \hline
 \multicolumn{5}{|c|}{\textbf{Probability of Quota Violations with $(\frac{q_1}{M}, \frac{q_2}{M}, \frac{q_3}{M}) \sim \mathrm{Dir}(2,2,2)$}}\\
 \hline
 
 Method & $M$ & Theoretical Probability & Sample Probability & $95\%$ Confidence Interval \\ \hline
 Huntington-Hill & $3$ & $0.13580$ & $0.13519$ & $( 0.13307, 0.13731)$ \\
 Adams & $3$ & $0.13580$ & $0.13486$ & $( 0.13274 , 0.13698 )$ \\
 Dean & $3$ & $0.13580$ & $0.13587$ & $(0.13375,0.13799)$ \\
 Jefferson & $3$ & $0.02670$ & $0.02627$ & $(0.02528,0.02726)$ \\
 Huntington-Hill & $5$ & $0.02276$ & $0.2327$ & $(0.02234,0.02420)$ \\
 Adams & $5$ & $0.03660 $ & $0.03729$ & $(0.03612,0.03846)$ \\
 Dean & $5$ & $0.02513$ & $0.02497$ & $(0.02400,0.02594)$ \\
 Jefferson & $5$ & $0.02889$ & $0.02862$ & $( 0.02759 , 0.02965 )$ \\
 Huntington-Hill & $10$ & $0.00386$ & $0.00391$ & $(0.00352,0.00430)$ \\
 Adams & $10$ & $0.03091$ & $0.0304$ & $(0.02934,0.03146)$ \\
 Dean & $10$ & $0.00537$ & $0.0052$ & $(0.00475,0.00565)$ \\
 Jefferson & $10$ & $0.03045$ & $0.03104$ & $(0.02997,0.03211)$ \\
 Huntington-Hill & $15$ & $0.00150$ & $0.00154$ & $(0.00130,0.00178)$ \\
 Adams & $15$ & $0.03107$ & $0.03195$ & $(0.03086,0.03304)$ \\
 Dean & $15$ & $0.00227$ & $0.00242$ & $(0.00212,0.00272)$ \\
 Jefferson & $15$ & $0.03094$ & $0.03015$ & $(0.02909,0.03121)$ \\

 \hline
\end{tabular}
\end{center}

\begin{remark}
As a Dirichlet $(1,1,1)$ distributions is a uniform distribution on the simplex $\{(x_1, x_2, x_3) \mid x_1 + x_2 + x_3 = 1 \}$, the case $(\frac{q_1}{M}, \frac{q_2}{M}, \frac{q_3}{M}) \sim \mathrm{Dir}(1,1,1)$ is equivalent to the case presented in Sections \ref{sample1section} and \ref{sample2section}.
 
\end{remark}

\begin{center}
\begin{tabular}{|c|c|c|c|c|}
 \hline
 \multicolumn{5}{|c|}{\textbf{Probability of Quota Violations with $p_i$ IID and $p_i \sim \exp (1)$}}\\
 \hline
 
 Method & $M$ & Theoretical Probability & Sample Probability & $95\%$ Confidence Interval \\ \hline
 Huntington-Hill & $3$ & $0.33333$ & $0.33413$ & $(0.33121,0.33705)$ \\
 Adams & $3$ & $0.33333$ & $0.33298$ & $(0.33006,0.33590)$ \\
 Dean & $3$ & $0.33333$ & $0.33379$ & $(0.33087,0.33671)$ \\
 Jefferson & $3$ & $0.02222$ & $0.02199$ & $(0.02108,0.02290)$ \\
 Huntington-Hill & $5$ & $0.13092$ & $0.13225$ & $(0.13015,0.13434)$ \\
 Adams & $5$ & $0.16000$ & $0.15939$ & $(0.15712,0.16166)$ \\
 Dean & $5$ & $0.13714$ & $0.13545$ & $(0.13333,0.13757)$ \\
 Jefferson & $5$ & $0.03943 $ & $0.03942$ & $(0.03821,0.04063)$ \\
 Huntington-Hill & $10$ & $0.04904$ & $0.04949$ & $(0.04815,0.05083)$ \\
 Adams & $10$ & $0.11057$ & $0.11116$ & $(0.10921,0.11311)$ \\
 Dean & $10$ & $0.05581$ & $0.05663$ & $(0.05520,0.05806)$ \\
 Jefferson & $10$ & $0.05671$ & $0.05642$ & $(0.05499,0.05785)$ \\

 \hline
\end{tabular}
\end{center}

\begin{remark}
The probabilities in the case $p_i$ are independent and identically distributed random variables with $p_i \sim \exp(\lambda)$ are the same as the probabilities presented in sections \ref{sample1section} and \ref{sample2section}. This may be shown by computing the associated probability density function for $(q_1, q_2)$.
 
\end{remark}

\begin{center}
\begin{tabular}{|c|c|c|c|c|}
 \hline
 \multicolumn{5}{|c|}{\textbf{Probability of Quota Violations with $p_i$ IID and $p_i \sim U (0,1000)$}}\\
 \hline
 
 Method & $M$ & Theoretical Probability & Sample Probability & $95\%$ Confidence Interval \\ \hline
 Huntington-Hill & $3$ & $0.125$ & $0.12422 $ & $(0.122218, 0.12626)$ \\
 Adams & $3$ & $0.125$ & $0.12399$ & $(0.12195,0.12603)$ \\
 Dean & $3$ & $0.125$ & $0.1251$ & $(0.12305,0.12715)$ \\
 Jefferson & $3$ & $0.01389$ & $0.01341$ & $(0.01270,0.1412)$ \\
 Huntington-Hill & $5$ & $0.03923$ & $0.03966 $ & $(0.03845, 0.04087)$ \\
 Adams & $5$ & $0.05903$ & $0.05944$ & $(0.05797,0.06091)$ \\
 Dean & $5$ & $0.04360$ & $0.04444$ & $(0.04316,0.04572)$ \\
 Jefferson & $5$ & $0.01819$ & $0.01802$ & $(0.01720,0.01884)$ \\
 Huntington-Hill & $10$ & $0.01363$ & $0.01357$ & $(0.01285, 0.01429)$ \\
 Adams & $10$ & $0.04015$ & $0.03966$ & $(0.03845,0.04087)$ \\
 Dean & $10$ & $0.01612$ & $0.01571$ & $(0.01494,0.01648)$ \\
 Jefferson & $10$ & $0.02308$ & $0.02278$ & $(0.02186,0.02370)$ \\
 Huntington-Hill & $15$ & $0.00802$ & $0.00807 $ & $(0.00752, 0.00862)$ \\
 Adams & $15$ & $0.03664$ & $0.03724$ & $(0.03607,0.03841)$ \\
 Dean & $15$ & $0.00953$ & $0.01006$ & $(0.00944,0.01068)$ \\
 Jefferson & $15$ & $0.02536$ & $0.02561$ & $(0.02463,0.02659)$ \\

 \hline
\end{tabular}
\end{center}


\end{document}